\documentclass[onecolumn, draftcls]{IEEEtran}

\usepackage{amsbsy}
\usepackage{amsmath}
\usepackage{mathtools}
\usepackage{amsfonts} 
\usepackage{amssymb}  
\usepackage{amsthm}   
\usepackage{amsxtra}
\usepackage{bm}
\usepackage{multirow}
\usepackage{nomencl}
\usepackage{color,colortbl}
\usepackage{dsfont}
\usepackage{enumerate}
\usepackage{fancybox}
\usepackage{fancyhdr}
\usepackage{graphics}
\usepackage{graphicx}
\usepackage{latexsym}
\usepackage{subfig}
\graphicspath{{./images/}}
\usepackage{algorithmic}

\definecolor{Dblue}{rgb}{0,0,1}
\definecolor{Dbrown}{rgb}{0.59,0.4,0}
\definecolor{Dred}{rgb}{0.64,0,0}
\definecolor{Dgreen}{rgb}{0,0.4,0}
\def\real{\mathbb R}
\def \bs {\boldsymbol}
\def \mr {\mathrm}
\def \tb {\textbf}

\def \gr {\nabla}

\newtheorem{theorem}{Theorem}[section]
\newtheorem{lemma}{Lemma}[section]
\newtheorem{definition}{Definition}[section]

\makenomenclature
\begin{document}

\title{Smoothing and Decomposition for Analysis Sparse Recovery}
\author{Zhao~Tan,~\IEEEmembership{Student Member,~IEEE},
        Yonina~C.~Eldar,~\IEEEmembership{Fellow,~IEEE},
        Amir~Beck,~
        and~Arye~Nehorai,~\IEEEmembership{Fellow,~IEEE}
\thanks{Z. Tan and A. Nehorai are with the Preston M. Green Department of Electrical and Systems Engineering Department, Washington University in St. Louis, St. Louis,
MO, 63130 USA. E-mail: \{tanz, nehorai\}@ese.wustl.edu.}
\thanks{Y. C. Eldar is with the Department of Electrical Engineering, Technion---Israel Institute of Technology, Haifa 32000, Isreal. E-mail: yonina@ee.technion.ac.il.}
\thanks{A. Beck is with the Department of Industrial Enigneering and Management, Technion---Israel Institute of Technology, Haifa 32000, Isreal. E-mail: becka@ee.technion.ac.il.}
\thanks{The work of Z. Tan and A. Nehorai was supported by the AFOSR Grant FA9550-11-1-0210,
NSF Grant CCF-1014908, and ONR Grant N000141310050. The work of Y. C. Eldar was supported in part by the Israel Science Foundation under Grant no.
170/10, in part by the Ollendorf Foundation, in part by a Magnet grant Metro450 from
the Israel Ministry of Industry and Trade, and in part by the Intel Collaborative Research Institute for Computational Intelligence (ICRI-CI). The work of A. Beck  was partially supported by the Israel Science Foundation under grant ISF No.253/12.}}
\maketitle

\begin{abstract}
We consider algorithms and recovery guarantees for the analysis sparse model in which the signal is sparse with respect to a highly coherent frame. We consider the use of a monotone version of  the fast iterative shrinkage-thresholding algorithm (MFISTA) to solve the analysis sparse recovery problem. Since the proximal operator in MFISTA does not have a closed-form solution for the analysis model, it cannot be applied directly. Instead,  we examine two alternatives based on smoothing and decomposition transformations that relax the original sparse recovery problem, and then implement MFISTA on the relaxed formulation. We refer to these two methods as smoothing-based and decomposition-based MFISTA. We analyze the convergence of both algorithms, and establish that smoothing-based MFISTA converges more rapidly when applied to general nonsmooth optimization problems. We then derive a performance bound on the reconstruction error using these techniques. The bound proves that our methods can recover a signal sparse in a redundant tight frame when the measurement matrix satisfies a properly adapted restricted isometry property. Numerical examples demonstrate the performance of our methods and show that smoothing-based MFISTA converges faster than the decomposition-based alternative in real applications, such as MRI image reconstruction.
\end{abstract}



\section{Introduction}
Low-dimensional signal recovery exploits the fact that many natural signals are inherently low dimensional, although they may have high ambient dimension. Prior information about the low-dimensional space can be exploited to aid in recovery of the signal of interest. Sparsity is one of the popular forms of prior information, and is the prior that underlies the growing field of compressive sensing \cite{CS1}-\nocite{CS2}\nocite{CS3}\cite{eldarbook}. Recovery of sparse inputs has found many applications in areas such as imaging, speech, radar signal processing, sub-Nyquist sampling and more. A typical sparse recovery problem has the following linear form:
\begin{equation}\label{eq:model}
\bs b=\bs A\bs x+\bs w,
\end{equation}
in which $\bs A\in \mathbb{R}^{m\times n}$ is a measurement matrix, $\bs b \in \mathbb{R}^m$ is the measurement vector, and $\bs w \in \mathbb{R}^m$ represents the noise term. Our goal is to recover the signal $\bs x \in \mathbb{R}^n$. Normally we have $m < n$, which indicates that the inverse problem is ill-posed and has infinitely many solutions. To find a unique solution, prior information on $\bs x$ must be incorporated.

In the synthesis approach to sparse recovery, it is assumed that $\bs x$ can be expressed as a sparse combination of known dictionary elements, represented as columns of a matrix $\bs D\in \real^{n\times p}$ with $p \geq n$. That is $\bs x=\bs D\bs \alpha$ with $\bs \alpha$ sparse, i.e., the number of non-zero elements in $\bs \alpha$ is far less than the length of $\bs \alpha$. The main methods for solving this problem can be classified into two categories. One includes greedy methods, such as iterative hard thresholding \cite{IHTACS} and orthogonal matching pursuit \cite{OMP}. The other is based on relaxation-type methods, such as basis pursuit \cite{BP} and LASSO \cite{lasso}. These methods can stably recover a sparse signal $\bs \alpha$  when the matrix $\bs A\bs D$ satisfies the restricted isometry property (RIP) \cite{RIP}\nocite{RIPLasso}-\cite{TD10}. 

Recently, an alternative approach has became popular, which is known as the analysis method \cite{ACS1}, \cite{ACS2}. In this framework, we are given an analysis dictionary $\bs D^* (\bs D \in \mathbb{R}^{n\times p})$ under which $\bs D^* \bs x$ is sparse. Assuming, for example, that the $\ell_2$ norm of the noise $\bs w$ is bounded by $\varepsilon$, the recovery problem can be formulated as
\begin{equation}
\mathop{\min}_{\bs x \in \mathbb{R}^n}\|\bs D^* \bs x\|_0 \quad \text{subject to } \|\bs b-\bs A\bs x\|_2\leq\varepsilon.
\end{equation}
Since this problem is NP hard, several greedy algorithms have been proposed to approximate it, such as thresholding \cite{IHT} and  subspace pursuit \cite{GreedyACS}. 

Alternatively, the nonconvex $\ell_0$ norm can be approximated by the convex $\ell_1$ norm leading to the following relaxed problem, referred to as analysis basis pursuit (ABP):
\begin{equation}\label{ABP}
\quad \mathop{\min}_{\bs x \in \mathbb{R}^n}\|\bs D^* \bs x\|_1 \quad \text{subject to } \|\bs b-\bs A\bs x\|_2\leq\varepsilon.
\end{equation}
ABP is equivalent to the unconstrained optimization
 \begin{equation}\label{ALASSO}
 \quad \mathop{\min}_{\bs x \in \mathbb{R}^n}\frac{1}{2}\|\bs b-\bs A\bs x\|^2_2+\lambda \|\bs D^* \bs x\|_1,
\end{equation}
which we call analysis LASSO (ALASSO). The equivalence is in the sense that for any $\varepsilon>0$ there exists a $\lambda$ for which the optimal solutions of ABP and ALASSO are identical.

Both optimization problems ABP and ALASSO can be solved using interior point methods \cite{CVXopt}. However, when the problem dimension grows, these techniques become very slow since they require solutions of linear systems.
Another suggested approach is based on  alternating direction method of multipliers (ADMM) \cite{ADMM,ADMM_app}.  The efficiency of this method highly depends on nice structure of the matrices $\bs A$. Fast versions of first-order algorithms, such as the fast iterative shrinkage-thresholding algorithm (FISTA) \cite{FISTA}, are more favorable in dealing with large dimensional data since they do not require $\bs A$ to have any structure. The difficulty in directly applying first-order techniques to ABP \eqref{ABP} and ALASSO \eqref{ALASSO} is the fact that the nonsmooth term $\|\bs D^* \bs x\|_1$ is inseparable. A generalized iterative soft-thresholding algorithm was proposed in \cite{IC11} to tackle this difficulty. However, this approach converges relatively slow as we will show in one of our numerical examples. A common alternative is to transform the nondifferentiable problem into a smooth counterpart. In \cite{NESTA}, the authors used Nesterov's smoothing-based method \cite{Smoothing2} in conjunction with continuation (NESTA) to solve ABP \eqref{ABP}, under the assumption that the matrix $\bs A^*\bs A$ is an orthogonal projector. In \cite{MDJ07}, a smoothed version of ALASSO \eqref{ALASSO} is solved using a nonlinear conjugate gradient descent algorithm. To avoid imposing conditions on $\bs A$, we focus in this paper on the  ALASSO formulation \eqref{ALASSO}.

It was shown in \cite{Smoothing1} that one can apply any fast first-order method that achieves an $\varepsilon$-optimal solution within $O(\frac{1}{\sqrt{\varepsilon}})$ iterations, to an $\varepsilon$ smooth-approximation of the general nonsmooth problem and obtain an algorithm with $O(\frac{1}{\varepsilon})$ iterations. In this paper, we choose a monotone version of FISTA (MFISTA) \cite{MFISTA} as our fast first-order method, whose objective function values are guaranteed to be non-increasing. We apply the smoothing approach together with MFISTA leading to the smoothing-based MFISTA (SFISTA) algorithm. We also propose a decomposition-based MFISTA method (DFISTA) to solve the analysis sparse recovery problem. The decomposition idea is to introduce an auxiliary variable $\bs z$ in \eqref{ALASSO} so that MFISTA can be applied in a simple and explicit manner. This decomposition approach can be traced back to \cite{courant}, and has been widely used for solving total variation problems in the context of image reconstruction  \cite{decomp2}. 

Both smoothing and decomposition based algorithms for nonsmooth optimization problems are very popular in the literature. One of the main goals of this paper is to examine their respective performance. We show that SFISTA requires lower computational complexity to reach a predetermined accuracy. Our results can be applied to a general model, and are not restricted to the analysis sparse recovery problem. 

In the context of analysis sparse recovery, we show in Section~\ref{sdasr} that both smoothing and decomposition techniques solve the following optimization problem:
\begin{equation}\label{Ralasso}
\quad \mathop{\min}_{\bs x \in \real^n, \bs z \in \real^p}   \frac{1}{2} \bs \|\bs A\bs x-\bs b\|_2^2+\lambda \|\bs z\|_1+\frac{1}{2}\rho\|\bs z-\bs D^*\bs x\|_2^2, \\
\end{equation}
which we refer to as relaxed ALASSO (RALASSO). Another contribution of this paper is in proving recovery guarantees for  RALASSO \eqref{Ralasso}. With the introduction of the restricted isometry property adapted to $\bs D$ (D-RIP) \cite{ACS1}, previous work \cite{ACS1} \cite{ALASSO} studied recovery guarantees based on ABP \eqref{ABP} and ALASSO (\ref{ALASSO}). Here we combine the techniques in \cite{RIP} and \cite{ALASSO}, and obtain a performance bound on RALASSO \eqref{Ralasso}. We show that when $\sigma_{2s}<0.1907$ and $\|\bs D^*\bs A^*\bs w\|_\infty \leq \frac{\lambda}{2}$, the solution $\hat {\bs x}_\rho$ of RALASSO \eqref{Ralasso} satisfies
\begin{equation}\label{bd_RALASSO}
\|\hat {\bs x}_\rho-\bs x \|_2 \leq C_0\sqrt{s} \lambda+C_1\frac{\|\bs D^*\bs x-(\bs D^*\bs x)_s\|_1}{\sqrt{s}}+C_2 \frac{\lambda p}{\sqrt{s}\rho},
\end{equation}
where $p$ is the number of rows in $\bs D^*$, $C_0,C_1,C_2$ are constants, and we use $(\bs x)_s$ to denote the vector consisting of the largest $s$ entries of $|\bs x|$. As a special case, choosing $\rho \to \infty$ extends the bound in (\ref{bd_RALASSO})  and obtains the reconstruction bound for ALASSO \eqref{ALASSO} as long as $\sigma_{2s}<0.1907$, which improves upon the results of \cite{ALASSO}.

The paper is organized as follows. In Section \ref{smoothdecomp}, we introduce some mathematical preliminaries, and present SFISTA and DFISTA for solving RALASSO \eqref{Ralasso}. We analyze the convergence behavior of these two algorithms in Section~\ref{convergenceanalysis}, and show that SFISTA converges faster than DFISTA for a general model. Performance guarantees on RALASSO \eqref{Ralasso} are developed in Section~\ref{bound}. Finally, in Section \ref{numerical} we test our techniques on numerical experiments to demonstrate the effectiveness of our algorithms in solving the analysis recovery problem. We show that SFISTA performs favorably in comparison with DFISTA. A continuation method is also introduced to further accelerate the convergence speed.

Throughout the paper, we use capital italic bold letters to represent matrices and lowercase italic bold letters to represent vectors. For a given matrix $\bs D$, $\bs D^*$ denotes the conjugate matrix. We denote by $\bs D^*_{\mathcal{T}}$ the matrix that maintains the rows in $\bs D^*$ with indices in set $\mathcal{T}$, while setting all other rows to zero. Given a vector $\bs x$, $\|\bs x\|_1, \|\bs x\|_2$ are the $\ell_1,\ell_2$ norms respectively, $\|\bs x\|_0$  counts the number of nonzero components which will be referred to as the $\ell_0$ norm although it is not a norm, and $\|\bs x\|_\infty$ denotes the maximum absolute value of the elements in $\bs x$. We use $\bs x[i]$ to represent the $i$th element of $\bs x$. For a matrix $\bs A$, $\|\bs A\|_2$ is the induced spectral norm, and $\|\bs A\|_{p,q}=\mathop{\max} \frac{\|\bs A \bs x\|_p}{\|\bs x\|_q}.$ Finally, $\mr{Re}\langle \bs a, \bs b\rangle=\frac{\langle \bs a, \bs b\rangle +\langle \bs b,\bs a\rangle}{2}.$ We use $\mr{argmin} \{f(\bs x): \bs x=\bs z,\bs y\}$ to denote $\bs z$ or $\bs y$, whichever yields a smaller function value of $f(\bs x)$. 


\section{Smoothing and Decomposition for Analysis Sparse Recovery}\label{smoothdecomp}

In this section we present the smoothing-based and decomposition-based methods for solving ALASSO \eqref{ALASSO}. To do so, we first recall in Subsection \ref{subsec:prox-grad} some results related to proximal gradient methods that will be essential to our presentation and analysis.

\subsection{The Proximal Gradient Method}\label{subsec:prox-grad}
 We begin this section with the definition of Moreau's proximal (or ``prox") operator \cite{M65}, which is the key step in defining the proximal gradient method.

 Given a closed proper convex function $h:\mathbb{R}^n \rightarrow \mathcal{\mathbb{R}}\cup\{\infty\}$, the proximal operator of $h$ is defined by
\begin{equation}
\mr{prox}_h(\bs x)= \underset{\bs u\in \mathbb{R}^n}{\arg\min} \left \{h(\bs u)+\frac{1}{2}\|\bs u-\bs x\|^2_2\right \}.
\end{equation}
The proximal operator can be computed efficiently in many important instances. For example, it can be easily obtained when $h$ is an $l_p$ norm ($p \in [1,\infty))$, or an indicator of ``simple" closed convex sets such as the box, unit-simplex and the ball. More examples of proximal operators as well as a wealth of properties can be found, for example, in \cite{M62} \cite{BC11}.

The proximal operator can be used in order to compute smooth approximations of convex functions. Specifically, let $h$ be a closed, proper, convex function, and let $\mu>0$ be a given parameter. Define
\begin{equation}
h_{\mu}(\bs x) = \min_{\bs u \in \real^n} \left \{ h(\bs u)+\frac{1}{2 \mu} \|\bs u-\bs x\|_2^2 \right \}.
\end{equation}
It is easy to see that
\begin{equation}
h_{\mu}(\bs x) = h({\rm prox}_{\mu h}(\bs x))+\frac{1}{2 \mu} \|\bs x-{\rm prox}_{\mu h}(\bs x)\|_2^2.
\end{equation}
The function $h_{\mu}$ is called  the \textit{Moreau envelope of $h$} and has the following important properties (see \cite{M65} for further details):
\begin{itemize} \item $h_{\mu}(\bs x) \leq h(\bs x)$. \item $h_{\mu}$ is continuously differentiable and its gradient is Lipschitz continuous with constant $1/\mu$.
\item The gradient of $h_{\mu}$ is given by
\begin{equation}\label{proxgrad} \nabla h_{\mu}(\bs x) = \frac{1}{\mu} (\bs x- {\rm prox}_{\mu h}(\bs x)).\end{equation}
\end{itemize}

One important usage of the proximal operator is in the proximal gradient method that is aimed at solving the following composite problem:
\begin{equation} \label{general:model} \min_{\bs x \in \real^n} \{F(\bs x)+G(\bs x)\}.\end{equation}
Here $F: {\mathbb R}^n \rightarrow \real$ is a continuously differentiable convex function with a continuous gradient that has Lipschitz constant $L_{\nabla F}$:
$$ \|\nabla F (\bs x)-\nabla F(\bs y) \|_2 \leq L_{\nabla F} \|\bs x-\bs y\|_2, \quad \mbox{ for all } \bs x,\bs y \in {\mathbb R}^n,$$
\noindent  and $G: {\mathbb R}^n \rightarrow  {\mathbb R}\cup\{\infty\}$ is an extended-valued, proper, closed and convex function. The \textit{proximal gradient method} for solving (\ref{general:model}) takes the following form (see \cite{FISTA,BT09}):\\
\bigskip

\begin{tabular}{l}\hline
\\
 \textbf{Proximal Gradient Method For Solving (\ref{general:model})} \\
   \hline
   \\
    \tb{Input}: An upper bound $L\geq L_{\gr F}$. \\
    \tb{Step 0.} Take $\bs x_0 \in \real^n.$\\
    \tb{Step k.} ($k\geq 1$) \\
      \quad \quad \quad Compute $\bs x_k={\rm prox}_{\frac{1}{L} G} \left(\bs x_{k-1}-\frac{1}{L} \nabla F(\bs x_{k-1}) \right ).$\\
 \hline
 \\
\end{tabular}

The main disadvantage of the proximal gradient method is that it suffers from a relatively slow $O(1/k)$ rate of convergence of the function values. An accelerated version is the \textit{fast proximal gradient method}, also known in the literature as \textit{fast iterative shrinkage thresholding algorithm} (FISTA) \cite{FISTA,BT09}. When $G \equiv 0$, the problem is smooth, and FISTA coincides with Nesterov's optimal gradient method \cite{N83}. In this paper we implement a monotone version of FISTA (MFISTA) \cite{MFISTA}, which guarantees that the objective function value is non-increasing along the iterations.
\\
\bigskip

\begin{tabular}{l}\hline
\\
 \textbf{Monotone FISTA Method (MFISTA) For Solving (\ref{general:model})} \\
   \hline
   \\
    \tb{Input}: An upper bound $L\geq L_{\gr F}$. \\
    \tb{Step 0.} Take $\bs y_1=\bs x_0, t_1=1.$\\
    \tb{Step k.} ($k\geq 1$) Compute\\
  \quad \quad \quad$\bs z_k={\rm prox}_{\frac{1}{L}G} \left ( \bs y_k-\frac{1}{L} \nabla F(\bs y_k)\right ).$\\
  \quad \quad \quad$t_{k+1}=\frac{1+\sqrt{1+4t^2_k}}{2}$.\\
  \quad \quad \quad $\bs x_k=\mr {argmin}\{F(\bs x)+G(\bs x):\bs x=\bs z_k,\bs x_{k-1}\}$.\\
  \quad \quad \quad$\bs y_{k+1}=\bs x_k+\frac{t_k}{t_{k+1}}(\bs z_k-\bs x_k)+\frac{t_k-1}{t_{k+1}}(\bs x_k-\bs x_{k-1})$.\\
\\
 \hline
 \\
\end{tabular}

The rate of convergence  of the sequence generated by MFISTA is $O(1/k^2)$.
\begin{theorem}\label{fistaconverge}
\cite{MFISTA} Let $\{\bs x_k\}_{k \geq 0}$ be the sequence generated by MFISTA, and let $\hat{\bs x}$ be an optimal solution of (\ref{general:model}). Then
\begin{equation}
F(\bs x_k)+G(\bs x_k)-F(\hat{\bs x})-G(\hat{\bs x}) \leq \frac{2 L_{\nabla F} \|\bs x_0-\hat{\bs x}\|_2^2}{(k+1)^2}.
\end{equation}
\end{theorem}

\subsection{The General Nonsmooth Model}
The general optimization model we consider in this paper is
\begin{equation}\label{P0}
\quad \mathop{\min}_{\bs x\in \mathbb{R}^n} \{H(\bs x)=f(\bs x)+g(\bs D^* \bs x)\},
\end{equation}
where $f: \mathbb{R}^n \rightarrow \mathbb{R}$ is a  continuously differentiable convex function with a Lipschitz continuous gradient $L_{\gr f}$. The function $g: \mathbb{R}^p \rightarrow \mathbb{R} \cup \{ \infty\}$ is a closed, proper convex function which is not necessarily smooth, and $\bs D^* \in \mathbb{R}^{p\times n}$ is a given matrix. In addition, we assume that $g$ is Lipschitz continuous with parameter $L_g$:
\begin{equation*}
|g(\bs z)-g(\bs v)|\leq L_g\|\bs z-\bs v\|_2 \quad \mbox{ for all } \bs z, \bs v \in \mathbb{R}^p.
\end{equation*}
This is equivalent to saying that the subgradients of $g$ over $\mathbb{R}^p$ are bounded by $L_g$:
\begin{equation*}
\|g' (\bs z)\|_2\leq L_g \text{ for any } \bs x\in \mathbb{R}^n \text{ and } g'(\bs z)\in \partial g(\bs z).
\end{equation*}
An additional assumption we make throughout is that the proximal operator of $\alpha g(\bs z)$ for any $\alpha >0$ can be easily computed.

Directly applying MFISTA to \eqref{P0} requires computing the proximal operator of $g(\bs D^*\bs x)$. Despite the fact that we assume that it is easy to compute the  proximal operator of $g(\bs z)$, it is in general difficult to compute that of $\alpha g(\bs D^* \bs x)$. Therefore we need to transform the problem before utilizing MFISTA, in order to avoid this computation.

When considering ALASSO, $f(\bs x)=\frac{1}{2}\|\bs A\bs x-\bs b\|_2^2$ and $g(\bs D^*\bs x)=\lambda \|\bs D^*\bs x\|_1$. The Lipschitz constants are given by  $L_{\gr f}=\|\bs A\|_2^2$ and $L_g=\lambda \sqrt{p}$. The proximal operator of $\alpha g(\bs z)=\alpha \lambda \|\bs z\|_1$ can be computed as
\begin{eqnarray}
\label{prox_g}
\mr{prox}_{\alpha g}(\bs z) = \Gamma_{\lambda \alpha}(\bs z) = [|\bs z|-\lambda \alpha]_+\mr{sgn}(\bs z),
\end{eqnarray}
where for brevity, we denote the soft shrinkage operator by $\Gamma_{\lambda \alpha}(\bs z).$ Here $[\bs z]_+$ denotes the vector whose components are given by the maximum between $z_i$ and $0$. Note, however, that there is no explicit expression for the proximal operator of $g(\bs D^* \bs x)= \lambda \|\bs D^* \bs x\|_1$, i.e., there is no closed form solution to
\begin{equation}
\underset{\bs u\in \mathbb{R}^n}{\arg\min} \left \{\alpha\lambda\|\bs D^*\bs u\|_1+\frac{1}{2}\|\bs u-\bs x\|^2_2\right \}.
\end{equation}

In the next subsection, we introduce two popular approaches for transforming the problem \eqref{P0}: smoothing and decomposition. We will show in Sections \ref{smoothalg} and \ref{decompalg} that both transformations lead to algorithms which only require computation of the proximal operator of $g(\bs z)$, and not that of $g(\bs D^* \bs x)$.

\subsection{The Smoothing and Decomposition Transformations}\label{sdasr}
The first approach to transform \eqref{P0} is the smoothing method in which the nonsmooth function $g(\bs z)$ is replaced by its Moreau envelope $g_\mu(\bs z)$, which can be seen as a smooth approximation. By letting $\bs z=\bs D^* \bs x$ , the smoothed problem becomes
\begin{equation}\label{Pmu}
 \quad \mathop{\min}_{\bs x\in \mathbb{R}^n} \{H_\mu(\bs x)=f(\bs x)+g_\mu(\bs D^* \bs x)\},
\end{equation}
to which MFISTA can be applied since it only requires evaluating the proximal operator of $g(\bs z)$. From the general properties of the Moreau envelope, and from the fact that the norms of the subgradients of $g$ are bounded above by $L_g$, we can deduce that there exists some $\beta_1$, $\beta_2>0$ such that $\beta_1+\beta_2=L_g$ and $ g(\bs z)-\beta_1\mu \leq g_\mu(\bs z)\leq g(\bs z)+\beta_2 \mu \mbox{ for all } \bs z \in \mathbb{R}^p$ (see \cite{Smoothing1, Smoothing2}). This shows that a smaller $\mu$ leads to a finer approximation.

The second approach for transforming the problem is the decomposition method in which we consider:
\begin{equation}\label{Prho}
\mathop{\min}_{\bs x\in \mathbb{R}^n, \bs z\in \mathbb{R}^p} \left \{G_\rho(\bs x,\bs z)=f(\bs x)+g(\bs z)+\frac{\rho}{2}\|\bs z-\bs D^*\bs x\|^2_2 \right \}.
\end{equation}
With $\rho \to \infty$, this problem is equivalent to the following constrained formulation of the original problem \eqref{P0}:
\begin{align} 
&\min  \{f(\bs x)+g(\bs z)\} \notag \\ 
&\mbox{s.t.} \quad \bs z = \bs D^* \bs x,  \quad \bs x \in \real^n, \bs z \in \real^p.\ 
\end{align}

Evidently, there is a close relationship between the approximate models \eqref{Pmu} and \eqref{Prho}. Indeed, fixing $\bs x$ and minimizing the objective function of \eqref{Prho} with respect to $\bs z$ we obtain
\begin{align}
\min_{\bs x \in \real^n, \bs z\in \real^p} \left \{f(\bs x)+g(\bs z)+\frac{\rho}{2} \|\bs z -\bs D^* \bs x\|_2^2\right \} \notag \\
=\min_{\bs x \in \real^n} \left \{ f(\bs x)+g_{\frac{1}{\rho}} (\bs D^* \bs x)\right \}.
\end{align}
Therefore, the two models are equivalent in the sense that their optimal solution set (limited to $\bs x$) is the same when $\mu=\frac{1}{\rho}$. For analysis sparse recovery, both transformations lead to RALASSO \eqref{Ralasso}. However, as we shall see, the resulting smoothing-based and decomposition-based algorithms and their analysis are very different.

\subsection{The Smoothing-Based Method}\label{smoothalg}

 Since \eqref{Pmu} is a smooth problem we can apply an optimal first-order method such as MFISTA with $F=H_{\mu}=f(\bs x)+g_\mu(\bs D^*\bs x)$ and $G\equiv 0$ in equation (\ref{general:model}). The Lipschitz constant of $H_\mu$ is given by $L_{\gr f}+\frac{\|\bs D\|^2_2}{\mu}$, and according to (\ref{proxgrad}) the gradient of $\nabla g_{\mu}(\bs D^* \bs x)$ is equal to $\frac{1}{\mu} \bs D(\bs D^* \bs x -{\rm prox}_{\mu g} (\bs D^* \bs x))$. The expression $\mr{prox}_{\mu g}(\bs D^*\bs x)$ is calculated by first computing $\mr{prox}_{\mu g}(\bs z)$, and then letting $\bs z=\bs D^*\bs x$.

Returning to the analysis sparse recovery problem, after smoothing we obtain
\begin{equation}\label{smoothmain}
\min_{\bs x \in \real^n} \left \{ H_\mu(\bs x)=\frac{1}{2} \bs \|\bs A\bs x-\bs b\|_2^2+ g_\mu(\bs D^* \bs x) \right \},
\end{equation}
where
\begin{align*}
g_\mu(\bs D^*\bs x)=& \mathop{\min}_{\bs u} \left\{\lambda \|\bs u\|_1+\frac{1}{2\mu}\|\bs u-\bs D^*\bs x\|_2^2\right\}\\
=&\sum_{i=1}^p \lambda \mathcal H_{\lambda \mu}((\bs D^*\bs x)[i]).
\end{align*}
The function $\mathcal H_\alpha(x)$ with parameter $\alpha >0$ is the so-called Huber function \cite{H64}, and is given by
\begin{equation}
 \mathcal H_\alpha(x)=\left\{
  \begin{array}{ll}
   \frac{1}{2\alpha }x^2&  \text{if } |x| < \alpha \\
    |x|-\frac{\alpha}{2}&  \text{otherwise}.\\
  \end{array} \right.
\end{equation}
From (\ref{prox_g}), the gradient of $g_\mu (\bs D^*\bs x)$ is equal to
\begin{equation}
 \gr g_\mu(\bs D^*\bs x) =\frac{1}{\mu} \bs D(\bs D^* \bs x -\Gamma_{\lambda \mu}(\bs D^*\bs x)).
\end{equation}
Applying MFISTA to (\ref{smoothmain}), results in the SFISTA algorithm, summarized in Algorithm 1.

\bigskip

\begin{tabular}{l}\hline
\\
 \textbf{Algorithm1: Smoothing-based MFISTA (SFISTA)} \\
   \hline
   \\
    \tb{Input}: An upper bound $L\geq \|\bs A\|_2^2+\frac{\|\bs D\|^2_2}{\mu}$. \\
    \tb{Step 0.} Take $\bs y_1=\bs x_0, t_1=1.$\\
    \tb{Step k.} ($k\geq 1$) Compute\\
$\gr f(\bs y_k)=\bs A^*(\bs A\bs y_k-\bs b)$.\\
$\gr g_\mu(\bs D^*\bs x_{k-1}) =\frac{1}{\mu} \bs D(\bs D^* \bs x_{k-1} -\Gamma_{\lambda \mu}(\bs D^*\bs x_{k-1}))$.\\
$\bs z_k=\bs y_k-\frac{1}{L}(\gr f(\bs y_k)+\gr g_\mu(\bs D^*\bs x_{k-1}))$. \\
$t_{k+1}=\frac{1+\sqrt{1+4t^2_k}}{2}$.\\
$\bs x_k=\mr {argmin}\{H_\mu(\bs x):\bs x=\bs z_k,\bs x_{k-1}\}$.\\
$\bs y_{k+1}=\bs x_k+\frac{t_k}{t_{k+1}}(\bs z_k-\bs x_k)+\frac{t_k-1}{t_{k+1}}(\bs x_k-\bs x_{k-1})$.\\
\\
 \hline
 \\
\end{tabular}

\subsection{The Decomposition-Based Method}\label{decompalg}
We can also employ MFISTA on the decomposition model
\begin{equation}
\min _{\bs x \in \real^n, \bs z \in \real^p}  \{G_\rho(\bs x,\bs z)=F_\rho(\bs x,\bs z)+G(\bs x,\bs z)\},
\end{equation}
where we take the smooth part as $F_\rho(\bs x,\bs z)=f(\bs x)+\frac{\rho}{2}\|\bs z-\bs D^*\bs x\|^2_2$ and the nonsmooth part as $G(\bs x,\bs z)=g(\bs z)$. In order to apply MFISTA to \eqref{Prho}, we need to compute the proximal operator of $\alpha G$ for a given constant $\alpha>0$, which is given by
\begin{align}
\mr{prox}_{\alpha G}( \bs x, \bs z)
=\begin{pmatrix}
                   \bs x \\
                   {\rm prox}_{\alpha g}(\bs z) \\
                 \end{pmatrix}.
\end{align}

In RALASSO \eqref{Ralasso}, $G(\bs x, \bs z)=\lambda \|\bs z\|_1$ and $F_\rho(\bs x,\bs z)= \frac{1}{2} \bs \|\bs A\bs x-\bs b\|_2^2+\frac{1}{2}\rho\|\bs z-\bs D^*\bs x\|_2^2$. Therefore,
\begin{equation}
    \mr{prox}_{\alpha G}( \bs x, \bs z)=\begin{pmatrix}
                   \bs x \\
                   \Gamma_{\lambda \alpha}(\bs z) \\
                 \end{pmatrix}.
\end{equation}
The Lipschitz constant of $\gr F$ is equal to $(\|\bs A\|_2^2+\rho(1+\|\bs D\|_2^2))$. By applying MFISTA directly, we have the DFISTA algorithm, stated in Algorithm 2.
\bigskip

\begin{tabular}{l}\hline
\\
 \textbf{Algorithm 2:Decomposition-based MFISTA (DFISTA)} \\
   \hline
   \\
    \tb{Input}: An upper bound $L\geq (\|\bs A\|_2^2+\rho(1+\|\bs D\|_2^2))$. \\
    \tb{Step 0.} Take $\bs u_1=\bs x_0, \bs v_1=\bs z_0, t_1=1.$\\
    \tb{Step k.} ($k\geq 1$) Compute\\
$\gr_{\bs x}F_\rho(\bs u_k,\bs v_k)=\bs A^*(\bs A\bs u_k-\bs b)+\rho \bs D(\bs D^*\bs u_k-\bs v_k)$.\\
$\gr_{\bs z}F_\rho(\bs u_k,\bs v_k))=\rho (\bs v_k-\bs D^* \bs u_k)$.\\
$\bs p_k=\bs u_k-\frac{1}{L}\gr_{\bs x}F_\rho(\bs u_k,\bs v_k)$.
\\$\bs q_k= \Gamma_\frac{\lambda}{L}(\bs v_k-\frac{1}{L}\gr_{\bs z}F_\rho(\bs u_k,\bs v_k))$.\\
$t_{k+1}=\frac{1+\sqrt{1+4t^2_k}}{2}$.\\
$(\bs x_k,\bs z_k)$\\
$=\mr {argmin}\{G_\rho(\bs x,\bs z):(\bs x,\bs z)=(\bs p_k, \bs q_k),(\bs x_{k-1}, \bs z_{k-1})\}$.\\
$\bs u_{k+1}=\bs x_k+\frac{t_k}{t_{k+1}}(\bs p_k-\bs x_k)+\frac{t_k-1}{t_{k+1}}(\bs x_k-\bs x_{k-1})$.\\
$\bs v_{k+1}=\bs z_k+\frac{t_k}{t_{k+1}}(\bs q_k-\bs z_k)+\frac{t_k-1}{t_{k+1}}(\bs z_k-\bs z_{k-1})$.\\
\\
 \hline
 \\
\end{tabular}


\section{Convergence Analysis}\label{convergenceanalysis}
In this section we analyze the convergence behavior of both the smoothing-based  and decomposition-based methods. Convergence of  smoothing algorithms has been treated in \cite{Smoothing2, Smoothing1}. In order to make the paper self contained, we quote the main results here. We then analyze the convergence of the decomposition approach. Both methods require the same type of operations at each iteration: the computation of the gradient of the smooth function $f$, and of the proximal operator corresponding to $\alpha g$, which means that they have the same computational cost per iteration. However, we show that smoothing converges faster than decomposition based methods. Specifically, the smoothing-based algorithm is guaranteed to generate an $\varepsilon$-optimal solution within $O(1/\varepsilon)$ iterations, whereas the decomposition-based approach
 requires $O(1/\varepsilon^{1.5})$ iterations. We prove the results by analyzing SFISTA and DFISTA for the general problem \eqref{P0}, however, the same analysis can be easily extended to other optimal first-order methods, such as the one described in \cite{Smoothing2}.

\subsection{Convergence of the Smoothing-Based Method}
For SFISTA the sequence $\{\bs x_k\}$ satisfies the following relationship \cite{MFISTA}:
\begin{equation}\label{quad}
H_\mu(\bs x_k) - H_\mu(\hat {\bs x}_\mu) \leq \frac{2\left (L_{\gr f}+\frac{\|\bs D\|^2_2}{\mu} \right )\Lambda_1}{(k+1)^2},
\end{equation}
where $\Lambda_1$ is an upper bound on the expression $\|\hat {\bs x}_\mu-\bs x_0\|_2$ with $\hat {\bs x}_\mu$ being an arbitrary optimal solution of the smoothed problem \eqref{Pmu}, and $\bs x_0$ is the initial point of the algorithm. Of course, this rate of convergence is problematic since we are more interested in bounding the expression $H(\bs x_k)-\hat H$ rather than the expression $H_\mu(\bs x_k) - H_\mu(\hat {\bs x}_\mu)$, which is in terms of the smoothed problem. Here, $\hat H$ stands for the optimal value for original nonsmooth problem \eqref{P0}. For that, we can use the following result from \cite{Smoothing1}.

\begin{theorem}  \cite{Smoothing1} Let \{$\bs x_k$\} be the sequence generated by applying MFISTA to the problem \eqref{Pmu}. Let $\bs x_0$ be the initial point and let $\hat {\bs x}$ denote the optimal solution of  \eqref{P0}.  An $\varepsilon$-optimal solution of \eqref{P0}, i.e. $|H(\bs x_k)-H(\hat {\bs x})|\leq \varepsilon$, is obtained in the smoothing-based method using MFISTA after at most
\begin{equation}
 K=2\|\bs D\|_2 \sqrt{L_g \Lambda_1}\frac{1}{\varepsilon}+\sqrt{L_{\gr f} \Lambda_1} \frac{1}{\sqrt{\varepsilon}}
\end{equation}
iterations with $\mu$ chosen as
\begin{equation}
\mu=\sqrt{\frac{\|\bs D\|_2^2}{L_g}}\frac{\varepsilon}{\sqrt{\|\bs D\|_2^2 L_g}+\sqrt{\|\bs D\|_2^2 L_g +L_{\gr f}\varepsilon}},
\end{equation}
in which $L_g$ and $L_{\gr f}$ are the Lipschitz constants of $g$ and the gradient function of $f$ in \eqref{P0}, and $\Lambda_1=\|\bs x_0-\hat {\bs x}_\mu\|_2$. We use $\hat {\bs x}_\mu$ to denote the optimal solution of problem \eqref{Pmu}.
\end{theorem}

\noindent \tb{Remarks:}
For analysis sparse recovery using SFISTA,  $L_g=\lambda p^{\frac{1}{2}}$ and $L_{\gr f}=\|\bs A\|_2^2$, which can be plugged into the expressions in the theorem.

\subsection{Convergence of the  Decomposition-Based Method}
 A key property of the decomposition model \eqref{Prho} is that its minimal value is bounded above by the optimal value $\hat H$ in the original problem \eqref{P0}.
\begin{lemma}\label{GHrelation}
Let $\hat G_{\rho}$ be the optimal value of problem \eqref{Prho} and $\hat H$ be the optimal value of problem \eqref{P0}. Then $\hat G_\rho \leq \hat H$.
\end{lemma}

\noindent \tb{Proof:}
The proof follows from adding the constraint $\bs z=\bs D^* \bs x$ to the optimization:
\begin{align}
\hat G_\rho= &\mathop{\min}_{\bs x\in \mathbb{R}^n, \bs z\in \mathbb{R}^p} \left\{f(\bs x)+g(\bs z)+\frac{\rho}{2}\|\bs z-\bs D^* \bs x\|^2_2\right\}\notag\\
\leq &\mathop{\min}_{\bs x\in \mathbb{R}^n, \bs z\in \mathbb{R}^p,\bs z=\bs D^*\bs x} \left\{f(\bs x)+g(\bs z)+\frac{\rho}{2}\|\bs z-\bs D^* \bs x\|^2_2\right\}\notag \\
=& \mathop{\min}_{\bs x\in \mathbb{R}^n} \left\{f(\bs x)+g(\bs D^*\bs x)\right\},
\end{align}
which is equal to $\hat H$.

The next theorem is our main convergence result establishing that an $\varepsilon$-optimal solution can be reached after $O(1/\varepsilon^{1.5})$ iterations. By assuming that the functions $f$ and $g$ are nonnegative, which is not an unusual assumption, we have the following theorem.
\begin{theorem} Let $\{\bs x_k, \bs z_k\}$ be the sequences generated by applying MFISTA to \eqref{Prho} with both $f$ and $g$ both being nonnegative functions. The initial point is taken as $(\bs x_0, \bs z_0)$ with $\bs z_0=\bs D^* \bs x_0$. Let $\hat {\bs x}$ denote the optimal solution of the original problem \eqref{P0}. An $\varepsilon$-optimal solution of problem \eqref{P0}, i.e. $|H(\bs x_k)-H(\hat {\bs x})|\leq \varepsilon$, is obtained using the decomposition-based method after at most
\begin{equation}
K = \mathrm{max} \left\{ \frac{16\sqrt{(1+\|\bs D\|^2\Lambda_2 H(\bs x_0))}L_g}{\varepsilon^{1.5}}, \frac{2\sqrt{L_{\gr f}\Lambda_2}}{\sqrt{\varepsilon}}\right\}
\end{equation}
iterations of MFISTA with $\rho$ chosen as
\begin{equation}
\rho=\left(\frac{L_g\sqrt{2H(\bs x_0)}K^2}{2(1+\|\bs D\|^2)\Lambda_2}\right)^{2/3}.
\end{equation}
Here $L_g$ and $L_{\gr f}$ are the Lipschitz constants for $g$ and the gradient function of $f$ in \eqref{P0}, and $\Lambda_2=\|\bs x_0-\hat {\bs x}_\rho\|^2_2+\|\bs z_0-\hat {\bs z}_\rho\|^2_2$. We use $\hat {\bs x}_\rho, \hat {\bs z}_\rho$ to denote the optimal solutions  to \eqref{Prho}.
\end{theorem}

\noindent \tb{Proof:} Since the monotone version of FISTA is applied we have
\begin{align}\label{mono}
&f(\bs x_k)+g(\bs z_k)+\frac{\rho}{2}\|\bs z_k-\bs D^* \bs x_k\|_2^2 \notag\\
=&G_\rho(\bs x_k,\bs z_k)\leq G_\rho(\bs x_0,\bs z_0)= f(\bs x_0)+g(\bs D^* \bs x_0)=H(\bs x_0).
\end{align}
With the assumption that $f$ and $g$ are nonnegative, it follows that
\begin{equation*}
\frac{\rho}{2}\|\bs z_k-\bs D^*\bs x_k\|^2_2\leq H(\bs x_0),
\end{equation*}
and therefore
\begin{equation}\label{zDxnorm}
\|\bs z_k-\bs D^*\bs x_k\|_2\leq \sqrt{\frac{2H(\bs x_0)}{\rho}}.
\end{equation}

The gradient of $f(\bs x)+\frac{\rho}{2}\|\bs z-\bs D^* \bs x\|_2^2$, is Lipschitz continuous with parameter $(L_{\gr f}+\rho(1+\|\bs D\|_2^2))$.  According to \cite{MFISTA}, by applying MFISTA, we obtain a sequence $\{(\bs x_k, \bs z_k)\}$ satisfying
\begin{equation*}
G_\rho(\bs x_k, \bs z_k)-\hat G_\rho \leq \frac{2(L_{\gr f}+\rho(1+\|\bs D\|^2_2))\Lambda_2}{k^2}.
\end{equation*}
Using lemma \ref{GHrelation} and the notation
\begin{align*}
A=2L_{\gr f}\Lambda_2, B=2(1+\|\bs D\|^2_2)\Lambda_2,
\end{align*}
we have
\begin{equation}\label{GF}
G_\rho(\bs x_k, \bs z_k)-\hat H\leq \frac{A+\rho B}{k^2}.
\end{equation}
We therefore conclude that
\begin{align*}
H(\bs x_k)=& f(\bs x_k)+g(\bs D^*\bs x_k)\\
                  =& f(\bs x_k)+g(\bs z_k)+g(\bs D^*\bs x_k)-g(\bs z_k)\\
                  \leq& G_\rho(\bs x_k,\bs z_k)+L_g\|\bs z_k-\bs D^* \bs x_k\|_2\\
                  \leq&\hat H+\frac{A+\rho B}{k^2}+L_g\|\bs z_k-\bs D^* \bs x_k\|_2\\
                  \leq& \hat H+\frac{A+\rho B}{k^2}+L_g\sqrt{\frac{2H(\bs x_0)}{\rho}}.
\end{align*}
The first inequality follows from the Lipschitz condition for the function $g$, the second inequality is obtained from (\ref{GF}), and the last inequality is a result of (\ref{zDxnorm}).

We now seek the ``best'' $\rho$ that minimizes the upper bound, or equivalently, minimizes the term
\begin{equation}
\frac{A+\rho B}{k^2}+L_g\sqrt{\frac{2H(\bs x_0)}{\rho}}=\frac{A}{k^2}+C\rho+\frac{D}{\sqrt{\rho}},
\end{equation}
where $C=\frac{B}{k^2}$ and $D=L_g\sqrt{2H(\bs x_0)}$. Setting the derivative to zero, the optimal value of $\rho$ is $\rho=\left (\frac{D}{2C} \right )^{2/3}$, and
\begin{equation}
H(\bs x_k) \leq \hat H+\frac{A}{k^2}+2C^{1/3}D^{2/3}.
\end{equation}
Therefore, to obtain an $\varepsilon$-optimal solution, it is enough that
\begin{equation}
\frac{A}{k^2}\leq \frac{\varepsilon}{2}, \quad \frac{2B^{1/3}D^{2/3}}{k^{2/3}}\leq \frac{\varepsilon}{2},
\end{equation}
or
\begin{align}
k\geq&\mr{max}\left\{\frac{4^{3/2}B^{1/2}D}{\varepsilon^{1.5}},\frac{\sqrt{2A}}{\sqrt{\varepsilon}}\right\}\notag \\
=& \mathrm{max} \left\{ \frac{16\sqrt{(1+\|\bs D\|^2\Lambda_2 H(\bs x_0))}L_g}{\varepsilon^{1.5}}, \frac{2\sqrt{L_{\gr f}\Lambda_2}}{\sqrt{\varepsilon}}\right\},
\end{align}
completing the proof.
\noindent

\tb{Remarks:}\\
\noindent 1. As in SFISTA, when treating the analysis sparse recovery problem, $L_g=\lambda p^{\frac{1}{2}}$ and $L_{\gr f}=\|\bs A\|_2^2$, which again can be plugged into the expressions in the theorem.\\
\noindent 2. MFISTA is applied in SFISTA and DFISTA to guarantee a mathematical rigorous proof, i.e. the existence of equation \eqref{mono}. In real application, FISTA without monotone operations can also be applied to yield corresponding smoothing and decomposition based algorithms.

Comparing the results of smoothing-based  and decomposition-based methods, we immediately conclude that the smoothing-based method is preferable. First, it requires only $O(1/\varepsilon)$ iterations to obtain an $\varepsilon$-optimal solution whereas the decomposition approach necessitates $O(1/\varepsilon^{3/2})$ iterations. Note that both bounds are better than the  bound $O(1/\varepsilon^2)$ corresponding to general sub-gradient schemes for nonsmooth optimization. Second, the bound in the smoothing approach depends on $\sqrt{L_g}$, and not on $L_g$, as when using   decomposition methods. This is important since, for example, when $g(\bs z)=\|\bs z\|_1$, we have $L_g=p^{\frac{1}{2}}$. In the smoothing approach the dependency on $p$ is of the form $p^{\frac{1}{4}}$ and not $p^{\frac{1}{2}}$, as when using the decomposition algorithm.


\section{ Performance Bounds}\label{bound}

We now turn to analyze the recovery performance of analysis LASSO when smoothing and decomposition are applied. As we have seen, both transformations lead to the same RALASSO problem in \eqref{Ralasso}. Our main result in this section shows that the reconstruction obtained by solving RALASSO is stable when $\bs D^*\bs x$ has rapidly decreasing coefficients and the noise in the model \eqref{eq:model} is small enough. Our performance bound also depends on the choice of parameter $\rho$ in the objective function. Before stating the main theorems, we first introduce a definition and some useful lemmas, whose proofs are detailed in the Appendix.

To ensure stable recovery, we require that the matrix ${\bs A}$ satisfies the D-RIP:
\begin{definition} (D-RIP) \cite{ACS1}.  The measurement matrix $\bs A$ obeys the restricted isometry property adapted to $\bs D$ with constant $\sigma_s$ if
\begin{equation}
(1-\sigma_s)\|\bs v\|_2^2 \leq \|\bs A\bs v\|_2^2 \leq (1+\sigma_s)\|\bs v\|_2^2
\end{equation}
holds for all $\bs v \in \Sigma_s=\{\bs y: \bs y=\bs D\bs x \text{ and } \|\bs x\|_0\leq s\}$. In other words, $\Sigma_s$ is the union of subspaces spanned by all subsets of $s$ columns of $\bs D$.
\end{definition}
The following lemma provides a useful inequality for matrices satisfying D-RIP.
\begin{lemma}\label{DRIP}
Let $\bs A$ satisfy the D-RIP with parameter $\sigma_{2s}$, and assume that $\bs u, \bs v \in \Sigma_s$. Then,
\begin{equation}
\mr {Re} \langle \bs A\bs u,\bs A\bs v\rangle \geq -\sigma_{2s} \|\bs u\|_2\|\bs v\|_2+ \mr{Re}\langle \bs u,\bs v\rangle.
\end{equation}
\end{lemma}

In the following, $\hat {\bs x}_\rho$ denotes the optimal solution of RALASSO \eqref{Ralasso} and $\bs x$ is the original signal in the linear model \eqref{eq:model}; we also use $\bs h$ to represent the reconstruction error $\bs h=\hat {\bs x}_\rho-\bs x$.  Let $\mathcal{T}$ be the indices of coefficients with $s$ largest magnitudes in the vector $\bs D^*\bs x$, and denote the complement of $\mathcal{T}$ by  $\mathcal{T}^c$.   Setting $\mathcal{T}_0=\mathcal{T}$, we decompose $\mathcal{T}_0^c$ into sets of size $s$ where $\mathcal{T}_1$ denotes the locations of the $s$ largest coefficients in $\bs D_{\mathcal{T}^c}^*\bs x$, $\mathcal{T}_2$ denote the next $s$ largest coefficients  and so on. Finally, we let $\mathcal{T}_{01}=\mathcal{T}_0\cup \mathcal{T}_1$.

Using the result of Lemma \ref{DRIP} and the inequality  $\|\bs D^*_{\mathcal{T}_{0}}\bs h\|_2+\|\bs D^*_{\mathcal{T}_{1}}\bs h\|_2 \leq \sqrt{2}\|\bs D^*_{\mathcal{T}_{01}}\bs h\|_2 $ since $\mathcal{T}_0$ and $\mathcal{T}_1$ are disjoint, we have the following lemma.
\begin{lemma}\label{DRIPproperty} \tb{(D-RIP property)}
Let $\bs h=\hat {\bs x}_\rho-\bs x$ be the reconstruction error in RALASSO \eqref{Ralasso}. We assume that $\bs A$ satisfies the D-RIP with parameter $\sigma_{2s}$ and $\bs D$ is a tight frame. Then,
\begin{align}
&\mr {Re} \langle \bs A\bs h,\bs A\bs D\bs D^*_{\mathcal{T}_{01}}\bs h \rangle \notag
\\ \geq&(1-\sigma_{2s}) \|\bs D^*_{\mathcal {T}_{01}}\bs h\|^2_2 -\sqrt{2}s^{-\frac{1}{2}}\sigma_{2s}\|\bs D^*_{\mathcal{T}_{01}}\bs h\|_2 \|\bs D^*_{\mathcal{T}^c}\bs h\|_1.
\end{align}
\end{lemma}

Finally, the lemmas below show that the reconstruction error $\bs h$ and $\|\bs D^*_{\mathcal{T}^c}\bs h\|_1$ can not be very large.

\begin{lemma}\label{optimalcon}\tb{(Optimality condition)}
The optimal solution $\hat {\bs x}_\rho$ for RALASSO \eqref{Ralasso} satisfies
\begin{align}
\|\bs D^*\bs A^*\bs A\bs h\|_{\infty} \leq \left (\frac{1}{2}+ \|\bs D^*\bs D\|_{1,1} \right )\lambda.
\end{align}
\end{lemma}

\begin{lemma}\label{cone} \tb{(Cone constraint)}
The optimal solution $\hat {\bs x}_\rho$ for RALASSO \eqref{Ralasso} satisfies the following cone constraint,
\begin{equation}
\|\bs D^*_{\mathcal{T}^c}\bs h\|_1 \leq \frac{\lambda}{\rho}p+3\|\bs D^*_{\mathcal{T}}\bs h\|_1+4\|\bs D^*_{\mathcal{T}^c}\bs x\|_1.
\end{equation}
\end{lemma}

We are now ready to state our main result.
\begin{theorem} \label{Main:thm1}
Let $\bs A$ be an $m \times n$ measurement matrix, $\bs D$ an arbitrary $n \times p$ tight frame, and let $\bs A$ satisfy the D-RIP with $\sigma_{2s} <0.1907$. Consider the measurement $\bs b=\bs A\bs x+\bs w$, where $\bs w$ is noise that satisfies  $\|\bs D^*\bs A^*\bs w\|_\infty \leq \frac{\lambda}{2}$. Then the solution $\hat {\bs x}_\rho$ to RALASSO \eqref{Ralasso} satisfies
\begin{equation}
\|\hat {\bs x}_\rho-\bs x \|_2 \leq C_0\sqrt{s} \lambda+C_1\frac{\|\bs D^*\bs x-(\bs D^*\bs x)_s\|_1}{\sqrt{s}}+C_2 \frac{\lambda p}{\sqrt{s}\rho},
\end{equation}
for the decomposition transformation and
\begin{equation}
\|\hat {\bs x}_\rho-\bs x \|_2 \leq C_0\sqrt{s} \lambda+C_1\frac{\|\bs D^*\bs x-(\bs D^*\bs x)_s\|_1}{\sqrt{s}}+C_2 \frac{\lambda\mu p}{\sqrt{s}},
\end{equation}
for the smoothing transformation. Here $(\bs D^*\bs x)_s$ is the vector consisting of the largest $s$ entries of $\bs D^*\bs x$ in magnitude, $C_1$ and $C_2$ are constants depending on $\sigma_{2s}$, and $C_0$ depends on $\sigma_{2s}$ and $\|\bs D^*\bs D\|_{1,1}$.
\end{theorem}

\noindent \tb{Proof: }
The proof follows mainly from the ideas in \cite{RIP}, \cite{ALASSO}, and proceeds in two steps. First, we try to show that $\bs D^*\bs h$ inside $\mathcal{T}_{01}$ is bounded by the terms of $\bs D^*\bs h$ outside the set $\mathcal{T}$. Then we show that $\bs D^*_{\mathcal{T}^c}\bs h$ is essentially small. 

From Lemma \ref{DRIPproperty},
\begin{align}\label{ADADDh1}
&\mr {Re} \langle \bs A\bs h,\bs A\bs D\bs D^*_{\mathcal{T}_{01}}\bs h \rangle \notag
\\ \geq& (1-\sigma_{2s})\|\bs D^*_{\mathcal{T}_{01}}\bs h\|^2_2-\sqrt{2} s^{-\frac{1}{2}}\sigma_{2s}\|\bs D^*_{\mathcal{T}_{01}}\bs h\|_2\|\bs D^*_{\mathcal{T}^c}\bs h\|_1.
\end{align}
Using the fact that $\mr{Re} \langle \bs x, \bs y \rangle \leq |\langle \bs x, \bs y \rangle| \leq \|\bs x\|_1\|\bs y\|_\infty$, we obtain that
\begin{align}\label{ADADDh2}
\mr {Re}\langle\bs  A\bs h,\bs A\bs D\bs D^*_{\mathcal{T}_{01}}\bs h\rangle =&\mr {Re} \langle \bs D^*\bs A^*\bs A\bs h,\bs D^*_{\mathcal{T}_{01}}\bs h \rangle \notag
\\ \leq& \|\bs D^*\bs A^*\bs A\bs h\|_{\infty}\|\bs D^*_{\mathcal{T}_{01}}\bs h\|_1\notag
\\ \leq& \sqrt{2s}c_0\lambda\|\bs D^*_{\mathcal{T}_{01}}\bs h\|_2,
\end{align}
with $c_0=\frac{1}{2}+ \|\bs D^*\bs D\|_{1,1}$. The second inequality is a result of Lemma \ref{optimalcon} and the fact that  $\|\bs D^*_{\mathcal{T}_{01}}\bs h\|_1\leq \sqrt{2s}\|\bs D^*_{\mathcal{T}_{01}}\bs h\|_2$, in which $2s$ is the number of nonzero terms in $\bs D^*_{\mathcal{T}_{01}}\bs h$. 
Combining (\ref{ADADDh1}) and (\ref{ADADDh2}), we get
\begin{equation}\label{Dh2}
\|\bs D^*_{\mathcal{T}_{01}}\bs h\|_2 \leq \frac{\sqrt{2s}\lambda c_0+\sqrt{2}s^{-\frac{1}{2}}\sigma_{2s}\|\bs D^*_{\mathcal{T}^c}\bs h\|_1}{1-\sigma_{2s}}.
\end{equation}
Then the second step bounds $\|\bs D^*_{\mathcal{T}^c}\bs h\|_1$. From (\ref{Dh2}),
\begin{align}\label{Dh3}
\|\bs D^*_{\mathcal{T}}\bs h\|_1 \leq& \sqrt{s} \|\bs D^*_{\mathcal{T}}\bs h\|_2 \leq \sqrt{s} \|\bs D^*_{\mathcal{T}_{01}}\bs h\|_2 \notag
\\ \leq&  \frac{\sqrt{2}\lambda sc_0+\sqrt{2}\sigma_{s}\|\bs D^*_{\mathcal{T}^c}\bs h\|_1}{1-\sigma_{2s}}.
\end{align}
Finally, using Lemma \ref{cone} and (\ref{Dh3}),
\begin{equation}
\|\bs D^*_{\mathcal{T}^c}\bs h\|_1 \leq \frac{\lambda}{\rho}p+\frac{3\sqrt{2}\lambda sc_0+3\sqrt{2}\sigma_{2s}\|\bs D^*_{\mathcal{T}^c}\bs h\|_1}{1-\sigma_{2s}}+4\|\bs D^*_{\mathcal{T}^c}\bs x\|_1.
\end{equation}
Since $\sigma_{2s}<0.1907$, we have $1-(1+3\sqrt{2})\sigma_{2s}>0$. Rearranging terms, the above inequality becomes
\begin{align}\label{Dh4}
&\|\bs D^*_{\mathcal{T}^c}\bs h\|_1 \notag
\\ \leq& \frac{1-\sigma_{2s}}{1-(1+3\sqrt{2})\sigma_{2s}}\frac{\lambda}{\rho}p+\frac{3\sqrt{2}\lambda sc_0+4(1-\sigma_{2s})\|\bs D^*_{\mathcal{T}^c}\bs x\|_1}{1-(1+3\sqrt{2})\sigma_{2s}}.
\end{align}

We now derive the bound on the reconstruction error. Using the results of (\ref{Dh2}) and (\ref{Dh4}), we get
\begin{align}
\|\bs h\|_2=&\|\bs D^*\bs h\|_2\leq \|\bs D^*_{\mathcal{T}_{01}}\bs h\|_2+\sum_{j\geq 2}\|\bs D^*_{\mathcal{T}_{j}}\bs h\|_2
\notag \\ \leq&  \frac{\sqrt{2s}\lambda c_0+\sqrt{2}s^{-\frac{1}{2}}\sigma_{2s}\|\bs D^*_{\mathcal{T}^c}\bs h\|_1}{1-\sigma_{2s}}+s^{-\frac{1}{2}}\|\bs D^*_{\mathcal{T}^c}\bs h\|_1
\notag \\=&\frac{c_0\lambda\sqrt{2s}}{1-\sigma_{2s}}+\frac{((\sqrt{2}-1)\sigma_{2s}+1)s^{-\frac{1}{2}}\|\bs D^*_{\mathcal{T}^c}\bs h\|_1}{1-\sigma_{2s}}
\notag \\ \leq& C_0\sqrt{s} \lambda+C_1\frac{\|\bs D^*\bs x-(\bs D^*\bs x)_s\|_1}{\sqrt{s}}+C_2 \frac{\lambda p}{\sqrt{s}\rho}.
\end{align}
The first equality follows from the assumption that $\bs D$ is a tight frame so that $\bs D\bs D^*=\bs I$. The first inequality is the result of the triangle inequality. The second inequality follows from \eqref{Dh2} and the fact that $\sum_{j\geq 2}\|\bs D^*_{\mathcal{T}_{j}}\bs h\|_2 \leq s^{-\frac{1}{2}}\|\bs D^*_{\mathcal{T}^c}\bs h\|_1$, which is proved in equation \eqref{boundtail} in the Appendix. The constants in the final result are given by
\begin{align*}
C_0=&\frac{4\sqrt{2}c_0}{1-(1+3\sqrt{2})\sigma_{2s}},
\\ C_1=&\frac{4((\sqrt{2}-1)\sigma_{2s}+1)}{1-(1+3\sqrt{2})\sigma_{2s}},
\\ C_2=&\frac{(\sqrt{2}-1)\sigma_{2s}+1}{1-(1+3\sqrt{2})\sigma_{2s}}. \quad
\end{align*}

To obtain the error bound for the smoothing transformation we replace $\rho$ with $1/\mu$ in the result. $\square$

Choosing $\rho \rightarrow \infty$ in RALASSO \eqref{Ralasso} leads to the ALASSO problem for which $\bs z=\bs D^*\bs x$. We then have the following result.
\begin{theorem} \label{Main:thm2}

Let $\bs A$ be an $m \times n$ measurement matrix, $\bs D$ an arbitrary $n \times p$ tight frame, and let $\bs A$ satisfy the D-RIP with $\sigma_{2s} <0.1907$. Consider the measurement $\bs b=\bs A\bs x+\bs w$, where $\bs w$ is noise that satisfies  $\|\bs D^*\bs A^*\bs w\|_\infty \leq \frac{\lambda}{2}$. Then the
 solution $\hat {\bs x}$ to ALASSO \eqref{ALASSO} satisfies
\begin{equation}
\|\hat {\bs x}-\bs x \|_2 \leq C_0\sqrt{s} \lambda+C_1\frac{\|\bs D^*\bs x-(\bs D^*\bs x)_s\|_1}{\sqrt{s}},
\end{equation}
where $(\bs D^*\bs x)_s$ is the vector consisting of the largest $s$ entries of $\bs D^*\bs x$ in magnitude, $C_1$ is a constant depending on $\sigma_{2s}$, and $C_0$ depends on $\sigma_{2s}$ and $\|\bs D^*\bs D\|_{1,1}$.
\end{theorem}

\noindent \tb{Remarks:}

\noindent 1. When the noise in the system is zero, we can set $\lambda$ as a positive value which is arbitrarily close to zero. The solution $\hat {\bs x}$ then satisfies $\|\hat {\bs x}-\bs x \| \leq C_1\frac{\|\bs D^*\bs x-(\bs D^*\bs x)_s\|_1}{\sqrt{s}}$, which parallels the result for the noiseless synthesis model in \cite{RIP}.

\noindent 2. When $\bs D^*$ is a tight frame, we have $\bs D\bs D^*=\bs I$. Therefore by letting $\bs v=\bs D^*\bs x$, we can reformulate the original analysis model as
\begin{equation}
\mathop{\min}_{\bs v}  \frac{1}{2} \bs \|\bs A \bs D\bs v-\bs b\|_2^2+\lambda \|\bs v\|_1.
\end{equation}
 Assuming that the noise term satisfies the $l_2$ norm constraint $\|\bs w\|_2 \leq \varepsilon$, we have 
\begin{align}
\|\bs D^*\bs A^*\bs w\|_{\infty}\leq \|\bs D^*\bs A^* \bs w\|_2\leq \|\bs D^* \bs A^*\|_2 \|\bs w\|_2\leq \varepsilon \|\bs D^* \bs A^*\|_2.
\end{align} 
When $\bs A$ satisfies D-RIP with $\sigma_{2s} < 0.1907$, by letting $\lambda=2 \varepsilon \|\bs D^* \bs A^*\|_2$ we have
\begin{equation}
\|\hat{\bs v}-\bs v \|_2 \leq \|\bs D^*\|_2 \|\hat {\bs x}-\bs x \|_2 \leq \tilde{C_0}\varepsilon+\tilde{C_1}\frac{\|\bs v-(\bs v)_s\|_1}{\sqrt{s}}.
\end{equation}
This result has a form similar to the reconstruction error bound shown in \cite{RIP}. However, the specific constants are different since in \cite{RIP} the matrix $\bs A\bs D$ is required to satisfy the RIP, whereas in our paper we require only that the D-RIP is satisfied.

\noindent 3. A similar performance bound is introduced in \cite{ALASSO} and shown to be valid when  $\sigma_{3s} <0.25$. Using Corollary 3.4 in \cite{sigma}, this is equivalent to  $\sigma_{2s} <0.0833$. Thus the results in Theorem \ref{Main:thm2} allow for a looser constraint on ALASSO recovery.

\noindent 4. The performance bound of Theorem \ref{Main:thm1} implies that a larger choice of $\rho$, or a smaller parameter $\mu$, leads to a smaller reconstruction error bound. This trend is intuitive since large $\rho$ or small $\mu$ results in smaller model inaccuracy. However, a larger $\rho$ or a smaller $\mu$ leads to a larger Lipschitz constant and thus results in slower convergence according to Theorem \ref{fistaconverge}. The idea of parameter continuation \cite{continuation} can be introduced to both $\rho$ and $\mu$ to accelerate the convergence while obtaining a desired reconstruction accuracy. More details will be given in the next section.


\section{Numerical Results}\label{numerical}
In the numerical examples, we use both randomly generated data and MRI image reconstruction to demonstrate that SFISTA performs better than DFISTA.  In the last example we also introduce a continuation technique to further speed up convergence of the smoothing-based method. We further compare SFISTA with the existing methods in \cite{ADMM_app, IC11,MDJ07} using MRI image reconstruction, and show its advantages.

\subsection{Randomly Generated Data in a Noiseless Case}

\begin{figure}[h]
  \centering
  \includegraphics[width=0.45\textwidth]{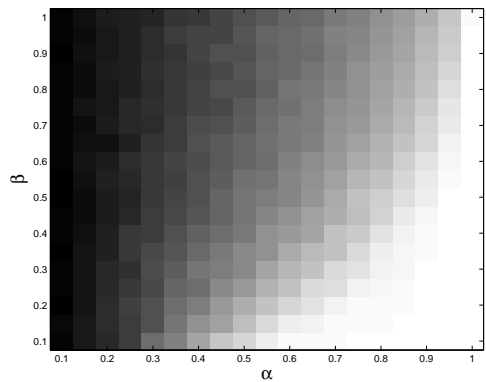}\\
  \caption{Reconstruction error of SFISTA}\label{smooth_pic}
\end{figure}

\begin{figure}[h]
  \centering
  \includegraphics[width=0.45\textwidth]{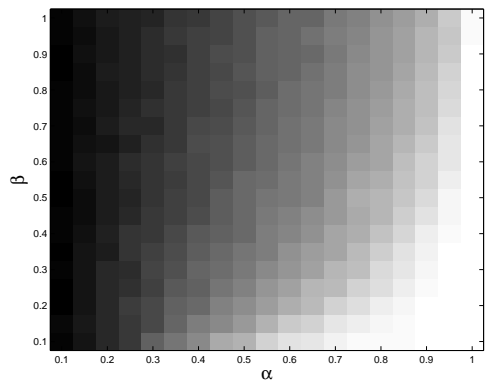}\\
  \caption{Reconstruction error of DFISTA}\label{decomp_pic}
\end{figure}

In this simulation, the entries in the $m \times n$ measurement matrix $\bs A$ were randomly generated according to a normal distribution. The $n \times p$ matrix $\bs D$ is a random tight frame. First we generated a $p \times n$ matrix whose elements follow an i.i.d Gaussian distribution. Then QR factorization was performed on this random matrix to yield the tight frame $\bs D$ with $\bs D\bs D^*=\bs I$ ($\bs D^*$ comprises the first $n$ columns from $\bs Q$, which was generated from the QR factorization). 

In the simulation we let $n = 120$ and $p = 144$, and we also set the values of $m$ and the number of zero terms named $l$ in $\bs D^*\bs x$ according to the following formula:
\begin{equation}
m = \alpha n, \quad l = n - \beta m.
\end{equation}
We varied $\alpha$ and $\beta$ from $0.1$ to $1$, with a step size $0.05$. We set $\lambda=0.004$, $\mu=10^{-3}\lambda^{-1}$ for the smoothing-based method, and $\rho=10^3\lambda$ for the decomposition-based method. For every combination of $\alpha$ and $\beta$, we ran a Monte Carlo simulation 50 times. Each algorithm ran for $3000$ iterations, and we computed the average reconstruction error. The reconstruction error is defined by $\frac{\|\hat{\bs x}-\bs x\|}{\|\bs x\|}$, in which $\hat {\bs x}$ is the reconstructed signal using smoothing or decomposition and $\bs x$ is the original signal in \eqref{eq:model}.

The average reconstruction error for smoothing and decomposition are plotted in Figs. \ref{smooth_pic} and \ref{decomp_pic}, respectively. White pixels present low reconstruction error whereas black pixels mean high error. Evidently, see that with same number of iterations, SFISTA results in a better reconstruction than DFISTA.

\subsection{MRI Image Reconstruction in a Noisy Case}
\begin{figure}[h]
  \centering
  \includegraphics[width=0.45\textwidth]{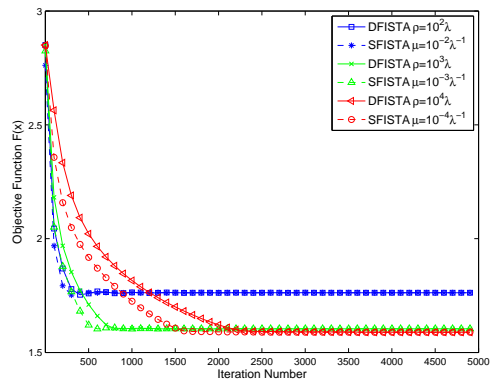}\\
  \caption{The objective function for MRI reconstruction on Shepp Logan. }\label{diffpara}
\end{figure}

\begin{figure}[h]
  \centering
  \includegraphics[width=0.45\textwidth]{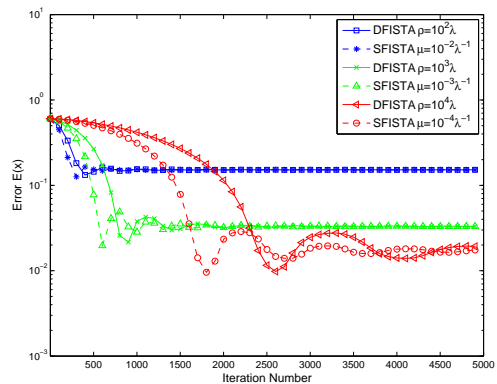}\\
  \caption{Reconstruction error for SFISTA and DFISTA with different parameters. }\label{semierror}
\end{figure}

The next numerical experiment was performed on a noisy $256 \times 256$ Shepp Logan phantom. The image scale was normalized to $[0,1]$. The additive noise followed a zero-mean Gaussian distribution with standard deviation $\sigma=0.001$. Due to the high cost of sampling in MRI, we only observed a limited number of radial lines of the phantom's 2D discrete Fourier transform. The matrix $\bs D^*$ consists of all vertical and horizontal gradients, which leads to a sparse $\bs D^* \bs x$. We let $\lambda=0.001$ in the optimization. We tested this MRI scenario with $\mu$ values of $10^{-2} \lambda^{-1}, 10^{-3} \lambda^{-1}, 10^{-4}\lambda^{-1}$ for SFISTA and  $\rho=10^2 \lambda, \rho=10^3 \lambda, 10^4 \lambda$ for DFISTA.  We took the samples along 15 radial lines to test these two methods.

In Fig.~\ref{diffpara} we plot the objective $\frac{1}{2} \bs \|\bs A\bs x-\bs b\|_2^2+\lambda \|\bs D^*\bs x\|_1$ as a function of the iteration number. It can be seen that the objective function of SFISTA decreases more rapidly than DFISTA. Furthermore, with smaller $\rho$ and larger $\mu$, DFISTA and SFISTA converge faster. Then we computed the reconstruction error. Here we see that smaller $\mu$ and larger $\rho$ lead to a more accurate reconstruction. We can see that SFISTA converges faster than DFISTA, which follows the convergence results in Section \ref{convergenceanalysis}.

\begin{figure}[h]
  \centering
  \includegraphics[width=0.45\textwidth]{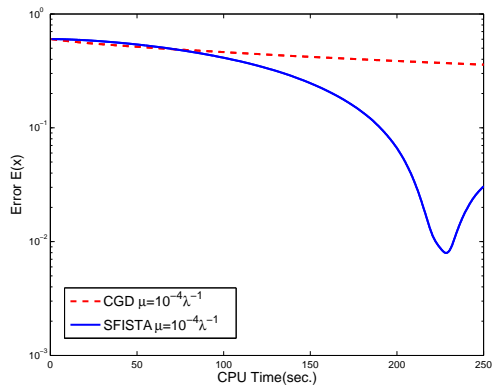}\\
  \caption{Reconstruction error for SFISTA and CGD with respect to CPU time. }\label{CGD}
\end{figure}

Next, we compared SFISTA with the nonlinear conjugate gradient descend (CGD) algorithm proposed in \cite{MDJ07}. The CGD also needs to introduce a smoothing transformation to approximate the term $\|\bs D^* \bs x\|_1$, and in this simulation the Moreau envelop with $\mu=10^{-4}\lambda^{-1}$ was used to smooth this term. We can see from Fig. \ref{CGD} that SFISTA converges faster than the CGD in terms of CPU time. CGD is slower because in each iteration, backtracking line-search is required, which reduces the algorithm efficiency.

\subsection{Acceleration by Continuation}

\bigskip

\begin{tabular}{l}\hline
\\
 \textbf{Algorithm 3: Continuation with SFISTA} \\
   \hline
   \\
    \tb{Input}: $\bs x$, the starting parameter $\mu=\mu_0$, \\
   \quad \quad\quad the ending parameter $\mu_f$ and $\gamma>1$. \\
    \tb{Step 1.} run SFISTA with $\mu$ and initial point $\bs x$. \\
    \tb{Step 2.} Get the solution $\bs x^*$ and let $\bs x=\bs x^*, \mu=\mu/\gamma$.\\
    \tb{Until.} \quad $\mu \leq \mu_f.$\\
\\
 \hline
 \\
\end{tabular}

To accelerate convergence and increase the accuracy of reconstruction, we consider continuation on the parameter $\mu$ for SFISTA, or on $\rho$ for DFISTA.
From Theorem \ref{Main:thm1}, we see that  smaller $\mu$ results in a smaller reconstruction error. At the same time, smaller $\mu$ leads to a larger Lipschitz constant $L_{\gr F}$ in Theorem \ref{fistaconverge}, and thus results in slower convergence. The idea of continuation is to solve a sequence of similar problems while using the previous solution as a warm start. Taking the smoothing-based method as an example, we can run SFISTA with $\mu_1\geq \mu_2 \geq \mu_3,\dots \geq \mu_f$. The continuation method is given in Algorithm 3. The algorithm for applying continuation on DFISTA is the same.

\begin{figure}[h]
  \centering
  \includegraphics[width=0.45\textwidth]{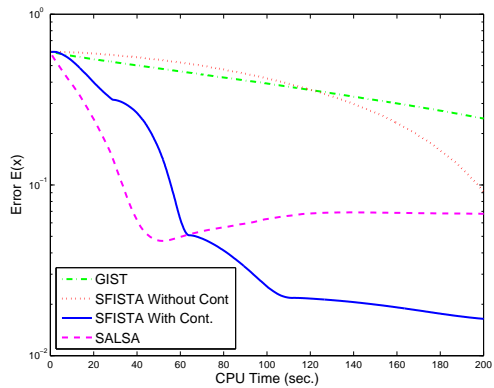}\\
  \caption{Convergence comparison among SFISTA with and without continuation, GIST and SALSA.}\label{recon_compare}
\end{figure}
\begin{figure}[h]
  \centering
  \includegraphics[width=0.45\textwidth]{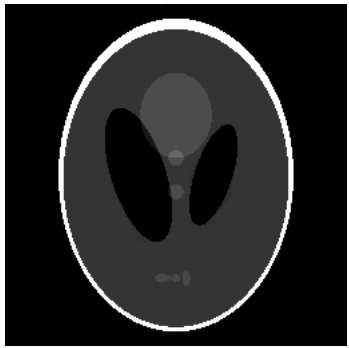}\\
  \caption{Reconstructed Shepp Logan with SFISTA using continuation.}\label{smoothcon}
\end{figure}

We tested the algorithm on the Shepp Logan image from the previous subsection with the same setting, using SFISTA with $\mu_f=10^{-4}\lambda^{-1}$ and standard SFISTA with $\mu=10^{-4}\lambda^{-1}$. We implemented the generalized iterative soft-thresholding algorithm (GIST) from \cite{IC11}. We also included an ADMM-based method, i.e. the split augmented Lagrangian shrinkage algorithm (SALSA) \cite{ADMM_app}. SALSA requires solving the proximal operator of $\|\bs D^*\bs x\|_1$, which is nontrivial. In this simulation, we implemented $40$ iterations of the Fast GP algorithm \cite{MFISTA} to approximate this proximal operator. Without solving the proximal operator exactly, the ADMM-based method can converge very fast while the accuracy of reconstruction is compromised as we show in Figure \ref{recon_compare}. In this figure we plot the reconstruction error for these four algorithms. It also shows that continuation helps speed up the convergence and exhibits better performance then GIST. The reconstructed Shepp Logan phantom using continuation is presented in Fig. \ref{smoothcon}, with reconstruction error  $3.17\%$. 


\section{conclusion}
In this paper, we proposed methods based on MFISTA to solve the analysis LASSO optimization problem.  Since the proximal operator in MFISTA for $\|\bs D^* \bs x\|_1$ does not have a closed-form solution, we presented two methods, SFISTA and DFISTA, using smoothing and decomposition respectively, to transform the original sparse recovery problem into a smooth counterpart. We analyzed the convergence of SFISTA and DFISTA and showed that SFISTA converges faster in general nonsmooth optimization problems. We also derived a bound on the performance for both approaches assuming a tight frame and D-RIP. Our methods were demonstrated via several simulations. With the application of parameter continuation, these two algorithms are suitable to solve large scale problems.

\appendix

\noindent \tb{Proof of Lemma \ref{DRIP}:}
Without loss of generality we assume that $\|\bs u\|_2=1$ and $\|\bs v\|_2=1$. By the definition of D-RIP, we have
\begin{align*}
\mr {Re} \langle \bs A\bs u,\bs A\bs v\rangle=&\frac{1}{4}\{\|\bs A\bs u+\bs A\bs v\|^2_2-\|\bs A\bs u-\bs A \bs v\|^2_2\}\\
\geq &\frac{1}{4}\{(1-\sigma_{2s})\|\bs u+\bs v\|^2_2-(1+\sigma_{2s})\|\bs u-\bs v\|^2_2\} \\
= &-\sigma_{2s}+\mr {Re}\langle \bs u,\bs v\rangle.
\end{align*}
Now it is easy to extend this equation to get the desired result.
\vspace{6pt}

\noindent \tb{Proof of Lemma \ref{DRIPproperty}:}
From the definition of $\mathcal{T}_j$ we have
$$\|\bs D^*_{\mathcal{T}_{j}}\bs h\|_2  \leq s^{-\frac{1}{2}}\|\bs D^*_{\mathcal{T}_{j-1}}\bs h\|_1$$
for all $j\geq 2$. Summing $j=2,3,\ldots$ leads to
\begin{equation}\label{boundtail}
\sum_{j\geq 2}\|\bs D^*_{\mathcal{T}_{j}}\bs h\|_2 \leq s^{-\frac{1}{2}} \sum_{j\geq 1}\|\bs D^*_{\mathcal{T}_{j}}\bs h\|_1 = s^{-\frac{1}{2}}\|\bs D^*_{\mathcal{T}^c}\bs h\|_1.
\end{equation}

Now, considering the fact that $\bs D$ is a tight frame, i.e., $\bs D\bs D^*=\bs I$, and that the D-RIP holds,
\begin{align*}
\mr {R}&\mr{e} \langle \bs A\bs h,\bs A\bs D\bs D^*_{\mathcal{T}_{01}}\bs h \rangle
\\=&\mr {Re} \langle \bs A\bs D \bs D^*_{\mathcal {T}_{01}}\bs h,\bs A\bs D\bs D^*_{\mathcal{T}_{01}}\bs h \rangle +\sum_{j\geq 2} \mr {Re} \langle \bs A\bs D \bs D^*_{\mathcal {T}_{j}}\bs h,\bs A\bs D\bs D^*_{\mathcal{T}_{01}}\bs h \rangle
\\ \geq& (1-\sigma_{2s}) \|\bs D \bs D^*_{\mathcal {T}_{01}}\bs h\|^2_2+\sum_{j\geq 2} \mr {Re} \langle \bs A\bs D \bs D^*_{\mathcal {T}_{j}}\bs h,\bs A\bs D\bs D^*_{\mathcal{T}_{0}}\bs h \rangle
\\ &+\sum_{j\geq 2} \mr {Re} \langle \bs A\bs D \bs D^*_{\mathcal {T}_{j}}\bs h,\bs A\bs D\bs D^*_{\mathcal{T}_{1}}\bs h \rangle
\end{align*}
Using the result from Lemma \ref{DRIP}, we can bound the last two terms in the above inequality; hence, we derive
\begin{align}\label{useful21}
\mr {R}&\mr{e} \langle \bs A\bs h, \bs A\bs D\bs D^*_{\mathcal{T}_{01}}\bs h \rangle \notag
\\ \geq&(1-\sigma_{2s}) \|\bs D \bs D^*_{\mathcal {T}_{01}}\bs h\|^2_2 + \sum_{j\geq 2} \mr {Re} \langle \bs D \bs D^*_{\mathcal {T}_{j}}\bs h,\bs D\bs D^*_{\mathcal{T}_{0}}\bs h \rangle \notag
\\&+\sum_{j\geq 2} \mr {Re} \langle \bs D \bs D^*_{\mathcal {T}_{j}}\bs h,\bs D\bs D^*_{\mathcal{T}_{1}}\bs h \rangle \notag
\\&-\sigma_{2s}\|\bs D\bs D^*_{\mathcal{T}_{0}}\bs h\|_2\sum_{j\geq 2}\|\bs D\bs D^*_{\mathcal{T}_{j}}\bs h\|_2 \notag
\\&-\sigma_{2s}\|\bs D\bs D^*_{\mathcal{T}_{1}}\bs h\|_2\sum_{j\geq 2}\|\bs D\bs D^*_{\mathcal{T}_{j}}\bs h\|_2 \notag
\\=&(1-\sigma_{2s}) \|\bs D \bs D^*_{\mathcal {T}_{01}}\bs h\|^2_2 + \mr {Re} \left \langle \sum_{j\geq 2} \bs D \bs D^*_{\mathcal {T}_{j}}\bs h,\bs D\bs D^*_{\mathcal{T}_{01}}\bs h \right \rangle \notag
\\ &-\sigma_{2s}(\|\bs D\bs D^*_{\mathcal{T}_{0}}\bs h\|_2+\|\bs D\bs D^*_{\mathcal{T}_{1}}\bs h\|_2)\sum_{j\geq 2}\|\bs D\bs D^*_{\mathcal{T}_{j}}\bs h\|_2
\end{align}
By definition of $\mathcal{T}_{j}$, we have
\begin{align*}
\mr{Re} \left \langle \sum_{j\geq 2} \bs D \bs D^*_{\mathcal {T}_{j}}\bs h,\bs D\bs D^*_{\mathcal{T}_{01}}\bs h \right \rangle&=\mr{Re}\langle \bs h-\bs D\bs D^*_{\mathcal{T}_{01}}\bs h,\bs D\bs D^*_{\mathcal{T}_{01}}\bs h \rangle
\\&= \|\bs D^*_{\mathcal{T}_{01}}\bs h\|^2_2- \|\bs D\bs D^*_{\mathcal{T}_{01}}\bs h\|^2_2.
\end{align*}
Combining this equation with (\ref{useful21}) results in
\begin{align*}
\mr {R}&\mr{e} \langle \bs A\bs h,\bs A\bs D\bs D^*_{\mathcal{T}_{01}}\bs h \rangle
\\ \geq&\|\bs D \bs D^*_{\mathcal {T}_{01}}\bs h\|^2_2-\sigma_{2s}\|\bs D \bs D^*_{\mathcal {T}_{01}}\bs h\|^2_2 + \|\bs D^*_{\mathcal{T}_{01}}\bs h\|^2_2- \|\bs D\bs D^*_{\mathcal{T}_{01}}\bs h\|^2_2
\\ &-\sigma_{2s}(\|\bs D\bs D^*_{\mathcal{T}_{0}}\bs h\|_2+\|\bs D\bs D^*_{\mathcal{T}_{1}}\bs h\|_2)\sum_{j\geq 2}\|\bs D\bs D^*_{\mathcal{T}_{j}}\bs h\|_2.
\end{align*}
Using the fact that when $\bs D$ is a tight frame, $\|\bs D \bs D^*_{\mathcal {T}_{01}}\bs h\|_2\leq \| \bs D^*_{\mathcal {T}_{01}}\bs h\|_2$, we have
\begin{align*}
\mr {R}&\mr{e}\langle \bs A\bs h,\bs A\bs D\bs D^*_{\mathcal{T}_{01}}\bs h \rangle
\\ \geq&(1-\sigma_{2s}) \|\bs D^*_{\mathcal {T}_{01}}\bs h\|^2_2 -\sigma_{2s}(\|\bs D^*_{\mathcal{T}_{0}}\bs h\|_2+\|\bs D^*_{\mathcal{T}_{1}}\bs h\|_2)\sum_{j\geq 2}\|\bs D^*_{\mathcal{T}_{j}}\bs h\|_2.
\end{align*}
Since $\|\bs D^*_{\mathcal{T}_{0}}\bs h\|_2+\|\bs D^*_{\mathcal{T}_{1}}\bs h\|_2 \leq \sqrt{2}\|\bs D^*_{\mathcal{T}_{01}}\bs h\|_2 $ (becuase  $\mathcal{T}_0$ and $\mathcal{T}_1$ are disjoint), we conclude that
\begin{align*}
\mr {R}&\mr{e} \langle \bs A\bs h,\bs A\bs D\bs D^*_{\mathcal{T}_{01}}\bs h \rangle
\\ \geq&(1-\sigma_{2s}) \|\bs D^*_{\mathcal {T}_{01}}\bs h\|^2_2 -\sqrt{2}\sigma_{2s}\|\bs D^*_{\mathcal{T}_{01}}\bs h\|_2\sum_{j\geq 2}\|\bs D^*_{\mathcal{T}_{j}}\bs h\|_2,
\end{align*}
which along with inequality (\ref{boundtail}) yields the desired result given by
\begin{align*}
\mr {R}&\mr{e} \langle \bs A\bs h,\bs A\bs D\bs D^*_{\mathcal{T}_{01}}\bs h \rangle
\\ \geq&(1-\sigma_{2s}) \|\bs D^*_{\mathcal {T}_{01}}\bs h\|^2_2 -\sqrt{2}s^{-\frac{1}{2}}\sigma_{2s}\|\bs D^*_{\mathcal{T}_{01}}\bs h\|_2 \|\bs D^*_{\mathcal{T}^c}\bs h\|_1.
\end{align*}

\vspace{6pt}
\noindent \tb{Proof of Lemma \ref{optimalcon}}: The subgradient optimality condition for RALASSO \eqref{Ralasso} can be stated as
\begin{equation}\label{sub1}
\bs A^*(\bs A\hat {\bs x}_\rho-\bs b)+\rho \bs D(\bs D^*\hat {\bs x}_\rho-\hat {\bs z}_\rho)=0,
\end{equation}
\begin{equation}\label{sub2}
\lambda \bs v + \rho (\hat {\bs z}_\rho-\bs D^*\hat {\bs x}_\rho)=0,
\end{equation}
where $\bs v$ is a subgradient  of the function $\|\bs z\|_1$ and consequently $\|\bs v\|_{\infty}\leq 1$. Combining (\ref{sub1}) and (\ref{sub2}), we have
\begin{equation*}
\bs A^*(\bs A\hat {\bs x}_\rho-\bs b)=\lambda\bs D \bs v.
\end{equation*}
Multiplying both sides by $\bs D^*$, we get
\begin{align}\label{DD11}
&\|\bs D^*\bs A^*(\bs A\hat{\bs x}_\rho-\bs b)\|_{\infty}\notag
\\&= \lambda \|\bs D^*\bs D \bs v\|_{\infty} \leq \lambda \|\bs D^*\bs D\|_{\infty,\infty}= \lambda \|\bs D^*\bs D\|_{1,1}.
\end{align}
The first inequality follows from the fact that $\|\bs v\|_{\infty}\leq 1$. With the assumption that $\|\bs D^*\bs A^*\bs w\|_\infty \leq \frac{\lambda}{2}$, and the triangle inequality, we have
\begin{align}\label{DAAh}
&\|\bs D^*\bs A^*\bs A\bs h\|_{\infty} \notag
\\ &\leq \|\bs D^*\bs A^*(\bs A\bs x-\bs b)\|_{\infty}+\|\bs D^*\bs A^*(\bs A\hat {\bs x}_\rho-\bs b)\|_{\infty}\notag
\\ &\leq \left (\frac{1}{2}+ \|\bs D^*\bs D\|_{1,1} \right )\lambda.
\end{align}
\vspace{6pt}

\noindent \tb{Proof of Lemma \ref{cone}}: Since $\hat {\bs x}_\rho$ and $\hat {\bs z}_\rho$ solve the optimization problem RALASSO \eqref{Ralasso}, we have,
\begin{align*}
&\frac{1}{2}\|\bs A\hat {\bs x}_\rho-\bs b\|^2_2+\lambda \|\hat {\bs z}_\rho\|_1+\frac{1}{2}\rho\|\bs D^*\hat {\bs x}_\rho-\hat {\bs z}_\rho\|^2_2
\\ &\leq \frac{1}{2}\|\bs A\bs x-\bs b\|_2^2+\lambda \|\bs D^*\bs x\|_1.
\end{align*}
Since $\bs b=\bs A\bs x+\bs w$ and $\bs h=\hat {\bs x}_\rho-\bs x$, it follows that
\begin{align*}
&\frac{1}{2}\|\bs A\bs h-\bs w\|_2^2+\lambda \|\hat {\bs z}_\rho\|_1+\frac{1}{2}\rho\|\bs D^*\hat {\bs x}_\rho-\hat {\bs z}_\rho\|^2_2
\\ &\leq \frac{1}{2}\|\bs w\|_2^2+\lambda \|\bs D^*\bs x\|_1.
\end{align*}
Expanding and rearranging the terms in the above equation, we get
\begin{align*}
&\frac{1}{2}\|\bs A\bs h\|_2^2+\lambda\|\hat {\bs z}_\rho\|_1+\frac{1}{2}\rho \|\bs D^*\hat {\bs x}_\rho-\hat {\bs z}_\rho\|_2^2
\\ &\leq \mr {Re} \langle \bs A\bs h,\bs w \rangle+\lambda \|\bs D^* \bs x\|_1,
\end{align*}
Using (\ref{sub2}) to replace the terms with $\hat {\bs z}_\rho$, we have
\begin{align*}
&\frac{1}{2}\|\bs A\bs h\|_2^2+\lambda \left \|\bs D^*\hat {\bs x}_\rho-\frac{\lambda}{\rho}\bs v\right \|_1+\frac{1}{2}\rho\left \|\frac{\lambda}{\rho}\bs v \right \|_2^2
\\&\leq \mr {Re} \langle \bs A\bs h,\bs w \rangle+\lambda\|\bs D^*\bs x\|_1.
\end{align*}
Since  $ \|\bs D^*\hat {\bs x}_\rho-\frac{\lambda}{\rho}\bs v\|_1 \geq \|\bs D^* \hat {\bs x}_\rho\|_1-\frac{\lambda}{\rho}\|\bs v\|_1$, we have
\begin{align}\label{AhDx}
&\frac{1}{2}\|\bs A\bs h\|_2^2+\lambda \|\bs D^*\hat {\bs x}_\rho\|_1\notag
\\ &\leq \frac{\lambda^2}{\rho} \|\bs v\|_1-\frac{\lambda^2}{2\rho}\|\bs v\|_2^2+\mr {Re} \langle \bs A\bs h,\bs w \rangle+\lambda\|\bs D^*\bs x\|_1\notag
\\ &\leq \frac{\lambda^2p}{2\rho}+\mr {Re}\langle \bs A\bs h,\bs w \rangle+\lambda\|\bs D^*\bs x\|_1.
\end{align}
The second inequality follows from the fact that $ \frac{\lambda^2}{\rho} \|\bs v\|_1-\frac{\lambda^2}{2\rho}\|\bs v\|_2^2$ is maximized when every element of $\bs v \in \real^p$ is $1$. Now, with the assumption that $\bs D$ is a tight frame, we have the following relation:
\begin{align*}
\mr {Re}\langle \bs A\bs h,\bs w\rangle+\lambda \|\bs D^* \bs x\|_1=&\mr {Re} \langle \bs D^*\bs h,\bs D^*\bs A^*\bs w\rangle+\lambda\|\bs D^*\bs x\|_1
\\ \leq& \|\bs D^* \bs h\|_1\|\bs D^*\bs A^*\bs w\|_\infty +\lambda \|\bs D^*\bs x\|_1.
\end{align*}
This inequality follows from the fact that $\mr{Re} \langle \bs x, \bs y \rangle \leq \|\langle \bs x, \bs y \rangle\| \leq \|\bs x\|_1\|\bs y\|_\infty$. Using the assumption that $\|\bs D^*\bs A^*\bs w\|_\infty \leq \frac{\lambda}{2}$, we get
\begin{equation}
 \mr {Re}\langle \bs A\bs h,\bs w\rangle+\lambda \|\bs D^* \bs x\|_1\leq \frac{\lambda}{2}\|\bs D^*\bs h\|_1+\lambda\|\bs D^*\bs x\|_1. \label{AwDx}
\end{equation}
Applying inequalities (\ref{AhDx}) and (\ref{AwDx}), we have
\begin{align*}
\lambda\|\bs D^*\hat {\bs x}_\rho\|_1\leq& \frac{1}{2}\|\bs A\bs h\|_2^2+\lambda \|\bs D^*\hat {\bs x}_\rho\|_1
\\ \leq& \frac{\lambda^2}{2\rho}p+\mr {Re}\langle \bs A\bs h,\bs w\rangle+\lambda\|\bs D^* \bs x\|_1
\\ \leq&  \frac{\lambda^2}{2\rho}p+\frac{\lambda}{2}\|\bs D^*\bs h\|_1+\lambda\|\bs D^*\bs x\|_1,
\end{align*}
which is the same as,
\begin{equation*}
\|\bs D^*\hat {\bs x}_\rho\|_1 \leq \frac{\lambda}{2\rho}p+\frac{1}{2}\|\bs D^*\bs h\|_1+\|\bs D^*\bs x\|_1.
\end{equation*}
Since we have $\bs h=\hat {\bs x}_\rho-\bs x$, it follows that
\begin{equation*}
\|\bs D^*\bs h+ \bs D^*\bs x\|_1 \leq \frac{\lambda}{2\rho}p+\frac{1}{2}\|\bs D^*\bs h\|_1+\|\bs D^*\bs x\|_1,
\end{equation*}
and hence
\begin{align*}
&\|\bs D^*_{\mathcal{T}}\bs h+ \bs D_{\mathcal{T}}^*\bs x\|_1+\|\bs D^*_{\mathcal{T}^c}\bs h+ \bs D_{\mathcal{T}^c}^*\bs x\|_1
\\ &\leq \frac{\lambda}{2\rho}p+\frac{1}{2}\|\bs D_{\mathcal{T}}^*\bs h\|_1+\frac{1}{2}\|\bs D_{\mathcal{T}^c}^*\bs h\|_1+\|\bs D_{\mathcal{T}}^*\bs x\|_1+\|\bs D_{\mathcal{T}^c}^*\bs x\|_1.
\end{align*}
Applying the triangle inequality to the left handside of above inequality, we results in
\begin{align*}
&-\|\bs D^*_{\mathcal{T}}\bs h\|_1+ \|\bs D_{\mathcal{T}}^*\bs x\|_1+\|\bs D^*_{\mathcal{T}^c}\bs h\|_1-\| \bs D_{\mathcal{T}^c}^*\bs x\|_1
\\ &\leq \frac{\lambda}{2\rho}p+\frac{1}{2}\|\bs D_{\mathcal{T}}^*\bs h\|_1+\frac{1}{2}\|\bs D_{\mathcal{T}^c}^*\bs h\|_1+\|\bs D_{\mathcal{T}}^*\bs x\|_1+\|\bs D_{\mathcal{T}^c}^*\bs x\|_1.
\end{align*}
After rearranging the terms, we have the following cone constraint,
\begin{equation}\label{Dhrelation}
\|\bs D^*_{\mathcal{T}^c}\bs h\|_1 \leq \frac{\lambda}{\rho}p+3\|\bs D^*_{\mathcal{T}}\bs h\|_1+4\|\bs D^*_{\mathcal{T}^c}\bs x\|_1.
\end{equation}


\bibliographystyle{IEEEtran}
\bibliography{IEEEabrv,acs}

\printnomenclature
\bibliographystyle{plain}

\end{document}